\documentclass[12pt]{amsart}
\usepackage{amsmath,amsfonts,amsthm,amscd,amssymb,graphicx}

\usepackage{amssymb}
\usepackage{hyperref}   
\usepackage{xcolor}

\def\rit{{\Bbb R}}

\def\nit{{\Bbb N}}

\def\tit{{\Bbb T}}

\def\eps{\varepsilon}

\def\cc{{\hbox{ c.c. }}}

\newtheorem{theorem}{Theorem}[section]

\newtheorem{e-proposition}[theorem]{Proposition}

\newtheorem{e-definition}[theorem]{Definition\rm}
\newtheorem{remark}{\it Remark\/}

\def\og{\leavevmode\raise.3ex\hbox{$\scriptscriptstyle\langle\!\langle$~}}
\def\fg{\leavevmode\raise.3ex\hbox{~$\!\scriptscriptstyle\,\rangle\!\rangle$}}

\def\beq{\begin{equation}}
\def\eeq{\end{equation}}

\title[]{ Boundary driven instabilities of Couette flows}

\author[D.Bian]{Dongfen Bian}
\address[D. Bian]{School of Mathematics and Statistics, Beijing Institute of Technology, Beijing 100081, China.}
\email{biandongfen@bit.edu.cn}

\author[E.Grenier]{Emmanuel Grenier}
\address[E. Grenier]{School of Mathematics and Statistics, Beijing Institute of Technology, Beijing 100081, China}
\email{Emmanuel.Grenier@ens-lyon.fr}

\author[N.Masmoudi]{Nader Masmoudi}
\address[N. Masmoudi]{NYUAD Research Institute, New York University Abu Dhabi, Saadiyat Island, Abu Dhabi, P.O. Box 129188, United Arab Emirates.\ 
Courant Institute of Mathematical Sciences, New York University, 251 Mercer Street New York, NY 10012 USA}
\email{masmoudi@cims.nyu.edu}
\author[W.Zhao]{Weiren Zhao}
\address[W. Zhao]{Department of Mathematics, New York University in Abu Dhabi, Saadiyat Island, P.O. Box 129188, Abu Dhabi, United Arab Emirates.}
\email{zjzjzwr@126.com, wz19@nyu.edu}

\begin{document}
\maketitle

%%%%%%%%%

%%%%%%%%%

\begin{abstract}
In this article, we prove that the threshold of instability of the classical Couette flow in $H^s$ for large $s$ is $\nu^{1/2}$.
The instability is completely driven by the boundary. The dynamic of the flow creates a Prandtl type boundary layer of width $\nu^{1/2}$ which is itself 
linearly unstable. This leads to a secondary instability which in turn  creates a sub-layer. 
\end{abstract}

%%%%%%%%%%%%%%%%%%%%%%%%%%%%%%%%%%%%%%%%%%%

\section{Introduction}

%%%%%%%%%%%%%%%%%%%%%%%%%%%%%%%%%%%%%%%%%%%

In this article, we consider the instability of the Couette flow with respect to the incompressible Navier-Stokes equations
\beq \label{NS1}
\partial_t u^\nu + (u^\nu \cdot \nabla) u^\nu - \nu \Delta u^\nu + \nabla p^\nu = f^\nu,
\eeq
\beq \label{NS2} 
\nabla \cdot u^\nu = 0,
\eeq
posed in the two dimensional infinite strip 
$$
\Omega_\rit = \{ (x,y) \quad | \quad x \in \rit, \quad -1  < y <  +1 \},
$$
or in the periodic strip
$$
\Omega_\tit = \{ (x,y) \quad | \quad x \in \tit, \quad -1  < y <  +1 \},
$$
 together with the initial condition
$u^\nu(0,\cdot,\cdot) = u_0^\nu$ and the boundary conditions 
\beq \label{NS3}
u^\nu(t,\cdot,-1) = (-1,0), \qquad u^\nu(t,\cdot,+1) = (1,0) .
\eeq
Let
$$
U_s(y) = (y,0)
$$
be the Couette flow.

When $\nu$ goes to $0$, we expect that boundary layers of size $\nu^{1/2}$, namely of Prandtl type,  appear near the lower and upper boundaries $y = \pm 1$.  As a consequence, even if we start with an initial data $u^\nu_0$ which has no boundary layer and whose derivatives are bounded, immediately, $u^\nu$ has large gradients near $y = \pm 1$. This leads to the following instability result in Sobolev spaces.

\begin{theorem} \label{theorem1}
Let $s$ be arbitrarily large and let $\beta > 0$. Then there exists a sequence of 
 solutions $u^\nu$ of \eqref{NS1}, \eqref{NS2}, and \eqref{NS3} with corresponding forces $f^\nu$ and a sequence of times
$T^\nu$, such that $u^\nu(0,\cdot,\pm 1)-(\pm 1,0) = 0$,
\beq \label{ii1}
\| u^\nu(0,\cdot,\cdot) - U_s(\cdot) \|_{H^s} \le \nu^{\beta},
\eeq
\beq \label{ii2}
 \| f^\nu \|_{L^\infty([0,T^\nu],H^s)} \le \nu^N,
\eeq
and such that
\beq \label{ii4}
\lim_{\nu \to 0}  \nu^{-\beta}\| \ \nabla u^\nu(T^\nu,\cdot,\cdot) \|_{H^1} = + \infty,
\eeq
\beq \label{ii5}
\lim_{\nu \to 0}  \nu^{-\beta} \| \nabla \times  u^\nu(T^\nu,\cdot,\cdot)  \|_{L^\infty} = + \infty,
\eeq
with $\lim_{\nu \to 0} T^\nu = 0$.
\end{theorem}

Note that this result is valid for any $\beta > 0$. Even for very small solutions (very large $\beta$), a boundary layer appears, leading to large amplifications of the gradient and of the vorticity of $u^\nu$.

In this article we study the instability of this boundary layer. To state an instability result, we have first to design a norm which is well-adapted to solutions having a boundary layer of size $\nu^{1/2}$, since, as is apparent in the previous theorem, Sobolev norms are not adequate.

A first possibility is to introduce the norms
\beq \label{norm1}
\|| u^\nu \||_s = \sum_{\alpha + \beta \le s}
\| \phi^{-\beta} \partial_x^\alpha \partial_y^\beta u^\nu \|_{L^\infty}
\eeq
where
$$
\phi(x,y) = 1 + {e^{- (y+1) / \delta \sqrt{\nu}}  \over \sqrt{\nu}} + {e^{- (1 - y) / \delta \sqrt{\nu}} \over \sqrt{\nu}},
$$
where $\delta > 0$ is some parameter. We can also replace $L^\infty$ by $L^2$ in (\ref{norm1}).
These norms appear to be very convenient to describe boundary layers.

However, in this article, we are interested in instability results, and an instability is  stronger the larger the space of arrival. Thus we will focus on a much larger norm
$$
\| u^\nu \|_s =  \sum_{\alpha + \beta \le s}
\sqrt{\nu}^\beta \|  \partial_x^\alpha \partial_y^\beta u^\nu \|_{L^\infty}.
$$
Functions with boundary layers of size $\sqrt{\nu}$ are uniformly bounded for these norms. This space also contains functions whose gradients are large "in the interior".
We will prove that the threshold of instability for these norms is $\beta = 1/2$.

We first prove an instability result when $\beta < 1/2$.
 The construction of the instability relies on a linear instability of the Prandtl boundary layer for the Rayleigh operator which is defined as follow.
For a given function $V_s(Y)$ and a given real number $\alpha\neq 0$, we define the Rayleigh operator to be
 $$
 Ray_{\alpha,V_s} = V_s - V_s'' (\partial_Y^2 - \alpha^2)^{-1},
 $$
where $ (\partial_Y^2 - \alpha^2)^{-1}$ is the inverse of $ (\partial_Y^2 - \alpha^2)$ with the Dirichlet boundary condition (see section \ref{preliminaries} for more details). 

Our instability result relies on the following claim.

\medskip

 {\it {\bf Claim:}
 Let 
 $$
\phi(t)  = - \cos { t } 
-   a \sinh {3  \over 2}  \cos {t \over 2}  -  a \sinh {1 \over 2}   \cos  {3  t \over 2},
$$
where 
$$
a = -\left(\sinh\frac{1}{2}+\sinh\frac{3}{2}\right)^{-1}.
$$
Let $w(t,Y)$ be the solution of
$$
\partial_t w - \partial_Y^2 w = 0
$$
with initial condition $w(0,\cdot) = 0$ and boundary conditions
$$
w(t,0) = \phi(t), \qquad
\hbox{and}\qquad  \lim\limits_{Y\to \infty}w(t,Y)=0.
$$
Then there exist $t_0 > 0$ and $\alpha_0 > 0$ such that the Rayleigh operator $Ray_{\alpha_0, V_{t_0}(Y)}$ associated with $\alpha_0$ and $V_{t_0}(Y)=\phi(t_0) + w(t_0,Y)$ has a non-real eigenvalue.}

\medskip

Note that $w(t,y)$ can be explicitly computed using double layers potentials. It is also possible to prove that $w(t,y)$ has an inflection point for some values of $t$, which implies that Rayleigh and Fjortoft criteria can not be applied.

Unfortunately there is no known method to prove that the Rayleigh operator has non real eigenvalues. However, it is numerically very simple and very classical to compute numerical approximations of the spectrum of Rayleigh operator. The problem reduces to the evaluation of the eigenvalues of a fixed matrix, or to the numerical integration of a simple non stiff ordinary differential equation. 

In this paper, using two different and independent approaches, we will give {\it numerical} evidences of the claim.

By accepting the claim, we prove the following instability result. 
\begin{theorem} \label{theomain} 
Assume that the claim is true. Let $s > 0$ and $N > 0$ be arbitrarily large, and let $\beta < 1/2$.
Then there exists a sequence of solutions $u^\nu$ of \eqref{NS1}, \eqref{NS2}, and \eqref{NS3} with corresponding forces $f^\nu$ and a sequence of times
$T^\nu$ such that
%such that $u^\nu(t,x,\pm 1)-(\pm 1,0) = 0$,
\beq \label{ii1b}
\| u^\nu(0,\cdot,\cdot) - U_s(\cdot) \|_{H^s} \le \nu^{\beta},
\eeq
\beq \label{ii2b}
 \| f^\nu \|_{C^1([0,T^\nu],H^s)} \le \nu^N,
\eeq
and such that
%\beq \label{ii3}
%\lim_{\nu \to 0} \| u^\nu(T^\nu,\cdot,\cdot) - U_s(\cdot) \|%_{H^{s'}} = + \infty,
%\eeq
\beq \label{ii4b}
\lim_{\nu \to 0} \nu^{1/2 - \beta} \| \nabla u^\nu(T^\nu,\cdot,\cdot) \|_{L^\infty} = + \infty,
\eeq
\beq \label{ii5b}
\lim_{\nu \to 0} \nu^{1/2 - \beta} \| \nabla \times  u^\nu(T^\nu,\cdot,\cdot) \|_{L^\infty} = + \infty,
\eeq
%\beq \label{ii6}
%\lim_{\nu \to 0}  \| \nabla \times  u^\nu(T^\nu,%\cdot,\cdot) \|_{L^\infty} = + \infty,
%\eeq
where $T^\nu \approx 1$.  At time $T^\nu$, $u^\nu$ has the scales $1$, $\nu^{1/2}$ and $\nu^{3/4 - \beta/2}$  in the $y$ variable.
\end{theorem}
\bigskip

In particular the vorticity is not bounded as $\nu$ goes to $0$. As physically expected, the instability is linked to the creation of a large vorticity near the boundary. More precisely, such an instability is due to the generation of a sub-layer with smaller size. Such a sub-layer does not appear in the case $\beta>1/2$, where we have the following theorem. 

\begin{theorem} \label{theo3} 
 Let us consider the Couette flow in the infinite strip $\Omega_\rit$ or the periodic strip $\Omega_\tit$.
Let $\beta > 1/2$. Let $s, N$ be large enough. Suppose that the initial data satisfies
$$
\|u^{\nu}(0,\cdot, \cdot)-U_s(\cdot)\|_{H^s}\leq \nu^{\beta}.
$$
%\begin{align*}
%    &\|u^{\nu}(0,\cdot, \cdot)-U_s(\cdot)\|_{H^s\leq \nu^{\beta},
%&\|f^{\nu}\|_{L^{\infty}H^{N}}\leq \nu^{N}.
%\end{align*}
Then, for any fixed $T > 0$, and any forcing satisfying 
\begin{align*}
    \| f^\nu \|_{C^1([0,T],H^s)} \le \nu^N,
\end{align*} 
it holds that
$$
\| u^\nu - U_s \|_{L^\infty([0,T],L^2)}
\le C \nu^{\beta},
$$
$$
\| u^\nu - U_s \|_{L^\infty([0,T],L^\infty)}
\le C \nu^{\beta},
$$
$$
\| \nabla (u^\nu-U_s) \|_{L^\infty([0,T],L^\infty)}
\le C \nu^{\beta - 1/2},
$$
for some constant $C$ independent on $\nu$.
In particular,
$$
\limsup_{\nu \to 0} \nu^{1/2-\beta}\| \nabla \times (u^\nu -U_s)\|_{L^\infty([0,T],L^\infty)} <\infty.
$$
%$$
%Then there exist $u^{app}$ and $T\lesssim 1$ such that for %$0\leq t<T$, the solutions $u^{\nu}$ to \eqref{NS1},% \eqref{NS2}, and \eqref{NS3} satisfy
%\begin{align*}
%    \|\nabla u^{app}\|_{L^{\infty}}\leq C\nu^{\beta-\frac{1}%{2}}
%\end{align*}
%and
%\begin{align*}
%    \|\nabla (u^{\nu}-u^{app})(t,\cdot, \cdot)\|_{L^{\infty}}\leq C\nu^{N-1}
%\end{align*}
%with $C$ independent of $\nu$. Here $u^{app}$ is obtained by %solving the linear system from the classical Prandtl Ansatz %\eqref{Ansatz}. 
\end{theorem}

With the same method, a similar result can be proved for the pipe flow $(0,0, 1 - (x^2+y^2))$ in the infinite tube
$$
{\mathcal T}_\rit = \{ (x,y,z) \quad | \quad z \in \rit, \quad |x|^2+y|^2  < 1 \Bigr\}
$$
or in the periodic tube
$$
{\mathcal T}_\tit = \{ (x,y,z) \quad | \quad z \in \tit, \quad |x|^2+|y|^2  < 1 \Bigr\}.
$$
Before going into the proof, let us discuss the physical mechanisms underlying these instabilities.
Let us assume that $u_0^\nu - U_s$ is of order $\nu^\beta$ for some $\beta > 0$.
When $\nu$ goes to $0$, we expect that $u^\nu$ is described by Prandtl Ansatz, namely that $u^\nu$ is of the form
$$
u^{\nu}(t,x,y) = u^{int}(t,x,y) + u^+ \Bigl( t, x, {y - 1 \over \sqrt{\nu}} \Bigr)
+ u^- \Bigl( t,x, {1 + y \over \sqrt{\nu}} \Bigr) + o(1)_{L^\infty},
$$
where $u^{int}$ describes the flow away from the boundary, and where $u^+$ and $u^-$, which rapidly decay in their last variable,
describe the boundary layers which appear at $y = \pm 1$. Moreover, we expect that $u^{int}$ is a solution of Euler equations,
whereas $u^\pm$ satisfy Prandtl equations.
In other words, even if we start with a perturbation which vanishes at $y = \pm 1$, we expect that the dynamics of the Navier-Stokes equations 
creates a boundary layer of width $\nu^{1/2}$ and amplitude $\nu^\beta$ by itself. The onset of this boundary layer leads to (\ref{ii4}) and (\ref{ii5}).

In physics, the stability of a  boundary layer is investigated through its  Reynolds number
$$
Re_{BL} = {U L \over \nu}
$$
where $U$ is the typical velocity in the boundary layer and $L$ its thickness. In our case,
$$
Re_{BL} = {\nu^\beta \nu^{1/2} \over \nu} = \nu^{\beta - 1/2}.
$$
If $Re_{BL} \to + \infty$, namely if $\beta < 1/2$,  we expect the boundary layer to be linearly and nonlinearly unstable.
On the contrary if $Re_{BL} \to 0$, namely if $\beta > 1/2$, these boundary layers are expected to be linearly stable, and more precisely, the $L^2$
norm of any perturbation in the boundary layer goes to $0$ and decreases monotonically.

We will prove that, if $\beta < 1/2$, this Prandtl boundary layer is indeed unstable, leading to a secondary instability and a sub-layer of size $\nu^{3/4 - \beta/2}$. This part relies on the numerical investigation of the claim previously stated.
 There exist arbitrarily small and smooth perturbations, of size $\nu^N$ in $H^s$ for
any arbitrarily large $N$ and $s$, which grow exponentially fast. These perturbations have two scales in $y$: the $\nu^{1/2}$ scale, which is the scale of the Prandtl
layer which creates them, but also a $\nu^{3/4-\beta/2}$ scale, which corresponds to a viscous sub-layer.
At these times, the solution $u^\nu$ has three scales, namely $1$ (scale of the domain), $\nu^{1/2}$ (scale of Prandtl boundary layer)
and $\nu^{3/4-\beta/2}$ (scale of the sub-layer created by the linear instability of Prandtl boundary layer).
%
%The stability of this sublayer is controlled by its Reynolds number $Re_{SL}$ defined by
%$$
%Re_{SL} = {U_{SL} L_{SL} \over \nu}
%$$
%where $U_{SL}$ is the typical velocity just outside the sublayer and $L_{SL}$ its width. It turns out that $U_{SL}$ is of order $\nu^\beta$,
%thus
%$$
%Re_{SL} \sim {\nu^\beta \nu^{3/4} \over \nu} \sim \nu^{\beta - 1/4}.
%$$
%Thus, if $1/4 < \beta < 1/2$, $Re_{SL}$ goes to $0$ as $\nu$ goes to $0$, and this sublayer is spectrally stable.
%If $\beta < 1/4$, this sublayer becomes itself unstable, leading to a third instability, which creates a new "sub-sublayer".
%

The stability of the Couette flow has been studied in the pioneer works of Kelvin \cite{Kelvin1887}, Rayleigh \cite{Rayleigh1880}, Orr \cite{Orr1907}, and Sommerfeld \cite{Sommerfeld1908}. It was suggested by Lord Kelvin \cite{Kelvin1887} that the stability/instability of the system is related to the size of the perturbation, and the threshold size is decreasing as $\nu\to 0$. With this perspective, the transition threshold problem, initially proposed by Trefethen et al. \cite{trefethen1993}, was later mathematically formulated by Bedrossian-Germain-Masmoudi \cite{BGM2017}:

{\it Given a norm $\|\cdot\|_X$, find a $\beta=\beta(X)$ so that
\begin{equation}\label{threshold}
  \begin{aligned}    
  &\|\omega_{in}\|_{X}\leq \nu^{\beta} \Rightarrow \text{stability},\\
&\|\omega_{in}\|_{X}\gg \nu^{\beta} \Rightarrow \text{instability}.    
  \end{aligned}
\end{equation}
}
Without boundary, for the 2D Couette flow, stability results have been established in different function spaces. For perturbations in Gevrey space (Gevrey-$\frac{1}{s}$, $\frac{1}{2}<s\leq 1$), $\beta\ge0$ indicates stability \cite{BMV2016}; for perturbations in Gevrey-$\frac{1}{s}$ with $0\leq s\leq \frac{1}{2}$, $\beta\geq \frac{1-2s}{3-3s}$ indicates stability \cite{li2022asymptotic}; for perturbations in Sobolev space ($H^\sigma$, $\sigma\ge2$), $\beta\ge \frac{1}{3}$ indicates stability \cite{MasmoudiZhao2019,weizhang2023}; and for perturbations in $L^2$ space, $\beta\ge \frac{1}{2}$ indicates stability \cite{MasmoudiZhao2020cpde}. In a recent work, Li-Masmoudi-Zhao \cite{LMZ2022} proved that the transition threshold for $L^2$ space is $\beta=\frac{1}{2}$. 

With boundary, two kinds of boundary conditions are studied. For the Navier-slip boundary condition, it was proved recently in \cite{bedrossian2024uniform} that for perturbations in Gevrey space (Gevrey-$\frac{1}{s}$, $\frac{1}{2}<s< 1$), $\beta\geq 0$ indicates stability. For the non-slip boundary condition, it was proved in \cite{chen2020transition} that for perturbations in $H^1$, $\beta\geq \frac{1}{2}$ indicates stability. Our result in the paper shows that $\beta=\frac{1}{2}$ is optimal threshold for Sobolev perturbations. 
\begin{remark}
   % The force $f^{\nu}$ in Theorem \ref{theomain} vanishes for $t\geq T^{\nu}$. 
    The initial perturbation $u^{\nu}(0)-U_s$ has compact support. With the same forcing and initial perturbations in the whole space setting, $\mathbb{T}\times \mathbb{R}$, it is easy to show that at $t=T^{\nu}$, $\|u^{\nu}(T^{\nu},\cdot, \cdot)-U_s\|_{H^s}\lesssim \nu^{\beta}$. Thus for $\beta\geq \frac{1}{3}$, the asymptotic stability holds. The instability in Theorem \ref{theomain} is boundary driven. 
\end{remark}

\section{Proof of Theorem \ref{theorem1} \label{boundary}}

%%%%%%%%%%%%%%%%%%%%%%%%%%%%%%%%%%%%%%%%%

We first prove Theorem \ref{theorem1} in the periodic case, namely in $\Omega_\tit$.

%%%%%%%%%%%%

\subsection{Preliminaries \label{preliminaries}}

%%%%%%%%%%%%

We first recall the classical Rayleigh and Orr-Sommerfeld equations which govern the spectral stability of shear flows for Euler and Navier-Stokes equations.

We take the Fourier transform in $x$ of these equations, with dual Fourier variable $\alpha$, and Laplace transform in time, with dual
variable $\lambda$. We define $c$ by $\lambda = - i \alpha c$.
The spectral stability of $U_s$ is linked to the existence of an non zero function $\psi$ and of a complex number $c$ with $\Im c > 0$,
solutions of the Orr-Sommerfeld equation
\beq \label{OS}
(U_s - c) (\partial_y^2 - \alpha^2) \psi - U_s'' \psi = \eps (\partial_y^2 - \alpha^2)^2 \psi,
\eeq
\beq \label{OS2}
\psi(\pm 1) = 0, \qquad \partial_y \psi(\pm 1) = 0,
\eeq
where $\eps = \nu / i \alpha$. When $\nu = 0$, these equations degenerate into Rayleigh equation
\beq \label{Ray}
(U_s - c) (\partial_y^2 - \alpha^2) \psi - U_s'' \psi = 0
\eeq
 with boundary condition
 \beq \label{Rayb}
 \psi(\pm 1) = 0.
 \eeq
According to Rayleigh's criterion, if there exists a solution $(c,\psi)$ of Rayleigh equation with $\Im c > 0$, then $U_s$ must have an inflexion point.
The reciprocal is however not true.

%%%%%%%

\subsection{Onset of the boundary layer}

%%%%%%%

We will consider initial data of the form
$$
u_0^\nu(x,y) = (y,0) + \nu^\beta u_1(x,y),
$$
where $u_1$ is a smooth vector filed which satisfies 
$$
u_1(x,-1) = u_1(x,+1) = 0
$$ 
for any $x$.
We  expect the solution $u^\nu$ to follow Prandtl Ansatz, namely to be of the form
\beq \label{Ansatz}
u^\nu(t,x,y) = (y,0) + \sum_{j = 1}^N \nu^{\beta_j} u_j(t,x,y) +  \nu^{\beta_j} u_j^{b}(t,x,Y) +  \nu^{\beta_j} u_j^{b,1}(t,x,Z)
\eeq
where
$$
Y = {y + 1 \over \nu^{1/2}}, \qquad Z = {y - 1 \over \nu^{1/2}} 
$$
and where $\beta_1 < \beta_2 < \cdots$ are various exponents of the form $k \beta + l/2$ with $(k,l) \in \nit^2$ and 
$(k,l) \ne (0,0)$.
In particular, $\beta_1 = \beta$.

In this Ansatz, the vector fields $u_j$ describe the behavior in the interior of the flow, whereas $u_j^b$ and $u_j^{b,1}$ are boundary layer correctors,
describing what happens close to the boundaries $y = -1$ and $y = 1$.

The equations on these various vector fields can be obtained by inserting (\ref{Ansatz}) into Navier-Stokes equations.
In particular, the leading profile $u_1(t,x,y)$ satisfies the linearized Euler equations
\beq \label{Euleru1}
\partial_t u_1 + (U_s \cdot\nabla) u_1 + (u_1 \cdot\nabla ) U_s + \nabla p_1 = 0,
\eeq
\beq \label{Euleru2}
\nabla \cdot u_1 = 0.
\eeq
Moreover, the various $u_j$ satisfy the same equation with some forcing term, coming from quadratic interactions between $u_l$ with $l < j$,
or from the diffusion of some $u_l$ with $l < j$.

 We now explicitly construct a solution of (\ref{Euleru1}) which does not vanish on the boundary.
 Let us fix some $\alpha > 0$.
Let us  look for $u_1$ under the form
$$
u_1(t,x,y) = \nabla^\perp( e^{i \alpha x} \psi_1(t,y)) + \cc
$$
where $\cc$ stands for ``complex conjugate".
Then the corresponding vorticity
$$
\omega_1(t,y) = (\partial_y^2 - \alpha^2) \psi_1(t,y)
$$
satisfies the transport equation
\beq \label{transportomega}
\partial_t \omega_1(t,y) + i \alpha y \omega_1(t,y) = 0.
\eeq
Hence
$$
\omega_1(t,y) = \omega_1(0,y) e^{ - i \alpha y t}  + \cc
$$
Let 
$$
G_\alpha(y',y)=\frac{-1}{\alpha\sinh \alpha}\left\{\begin{aligned}
    &\sinh\alpha(y-1)\sinh\alpha(y'+1)\quad -1\leq y'\leq y\leq 1\\
    &\sinh\alpha(y'-1)\sinh\alpha(y+1)\quad -1\leq y\leq y'\leq 1
\end{aligned}\right.
$$
be the Green function of $\partial_y^2 - \alpha^2$ with Dirichlet boundary conditions on $[-1,1]$.
Then
$$
\psi_1(t,y) = \int_{-1}^{+1} G_\alpha(z,y) \omega_1(0,z) e^{- i \alpha z t} \, dz.
$$
By construction of the Green function, 
$$
\psi_1(t,1) = \psi_1(t,-1) = 0
$$ 
for any $t \ge 0$. Moreover,
\beq \label{Dirac40}
\partial_y \psi_1(t,\pm 1) = \int_{-1}^{+1} \partial_y G_\alpha(z,\pm 1) \omega_1(0,z) e^{- i \alpha z t} \, dz.
\eeq
We take $\omega_1(0,y)$ of the form
\beq \label{Dirac}
\omega_1(0,y) = a_1 \delta_{b_1} + a_2 \delta_{b_2} + a_3 \delta_{b_3} 
\eeq
where $-1 < b_1 < b_2 < b_3 < 1$ and where $\delta_{b_i}$ is the Dirac mass at $b_i$. 
We then have
\beq \label{exprpsi}
\begin{aligned}
\partial_y \psi_1(t,\pm 1)
= &\partial_y G_\alpha(b_1,\pm 1) a_1 e^{-i \alpha b_1 t}
+ \partial_y G_\alpha(b_2,\pm 1) a_2 e^{- i \alpha b_2 t}\\
&+ \partial_y G_\alpha(b_3,\pm 1) a_3 e^{- i \alpha b_3 t}.
\end{aligned}
\eeq
At $t = 0$, $\partial_y \psi_1(0,\pm 1) = 0$  leads to the system
\beq \label{Dirac2}
\Bigl\{ \begin{array}{c} 
a_1 \sinh \alpha (b_1 + 1)  + a_2 \sinh \alpha (b_2 + 1)  + a_3 \sinh \alpha (b_3 + 1)  = 0, \cr
a_1 \sinh \alpha (b_1 - 1)  + a_2 \sinh \alpha (b_2 - 1)  + a_3 \sinh \alpha (b_3 - 1)  = 0. \cr
\end{array} 
\eeq
We choose 
$$
a_2 = {1 \over \sinh \alpha},  \qquad b_1 = - {1 \over 2}, \qquad b_2 = 0, \qquad b_3 =  {1 \over 2}.
$$ 
This leads to 
$$
a_1 \sinh {\alpha \over 2} + a_3 \sinh {3 \alpha \over 2} =  - 1,
$$
$$
a_1 \sinh {3 \alpha \over 2} + a_3 \sinh {\alpha \over 2} = - 1,
$$
which gives 
\beq \label{Dirac3}
a_1 = a_3 = -  \Bigl(  \sinh {\alpha \over 2}  + \sinh {3 \alpha \over 2} \Bigr)^{-1}. 
\eeq
Then, up to a factor $\sinh \alpha$,
$$ 
\partial_y \psi_1(t,-1) =  a_1 \sinh {3 \alpha \over 2} e^{ i {\alpha t / 2}}
+ a_2 \sinh \alpha + a_3 \sinh { \alpha \over 2} e^{-i \alpha t  / 2} . 
$$
Thus, $\partial_y \psi(t,\pm 1)$ is periodic, of period $4 \pi \alpha^{-1}$.
This leads to
$$
u_{1,h}(t,x,-1) = \Re \Bigl[  a_1 \sinh {3 \alpha \over 2} e^{ i {\alpha t / 2}} e^{i \alpha x} 
+ a_2 \sinh \alpha e^{i \alpha x} + a_3 \sinh { \alpha \over 2} e^{-i \alpha t  / 2} e^{i \alpha x} \Bigr],
$$
where $u_{1,h}(t,x,y)$ is the horizontal component of $u_1$.
In particular,
\beq \label{Dirac100}
u_{1,h}(t,x,-1) = 
\cos \alpha x  +  a_1 \sinh {3\alpha \over 2}  \cos \alpha \Bigl(  x + {t \over 2} \Bigr) +  a_3 \sinh { \alpha \over 2}   \cos \alpha \Bigl(  x - {t \over 2} \Bigr)  .
\eeq
In particular, $u_{1,h}(t,x,-1)$ does not identically vanish on the boundary. This implies that a boundary layer appears, namely that $u^b$ does not identically vanish. We study this boundary layer in the next section.

By continuity we can choose a smooth initial vorticity $\omega_1(0,y)$ such that $u_{1,h}(t,x,-1)$ remains arbitrarily close to (\ref{Dirac100}).
More precisely, let $\chi(y)$ be a smooth non-negative function, supported in $[-1,1]$ with unit integral. We replace (\ref{Dirac})
by
\beq \label{Diracbis}
\omega_1(0,y,\mu) = {a_1(\mu) \over \mu} \chi \Bigl( {y - b_1 \over \mu} \Bigr) 
+ {a_2(\mu) \over  \mu} \chi \Bigl( {y - b_2 \over \mu} \Bigr)  
+ {a_3(\mu) \over \mu} \chi \Bigl( {y - b_3 \over \mu} \Bigr)  
\eeq
where $\mu > 0$ is some small parameter.
Let $\psi_1(0,y,\mu)$ be the corresponding stream function.
We have
\beq \label{Dirac35}
\partial_y \psi_1(t,\pm 1,\mu) = \int_{-1}^{+1} \partial_y G_\alpha(z,\pm 1) \omega_1(0,z,\mu) e^{- i \alpha z t} \, dz.
\eeq
Then $\partial_y \psi_1(0,\pm 1,\mu) = 0$ leads to
\beq \label{Dirac10}
\sum_{k=1}^3 a_k(\mu) A_k^\pm(\mu) = 0
\eeq
where
$$
A_k^\pm(\mu) = \mu^{-1} \int_{-1}^{+1}\chi \Bigl( {z - b_k \over \mu} \Bigr) 
\partial_y G_\alpha(z, \pm 1) \, dz.
$$
As $\mu$ goes to $0$, the coefficients $A_k^\pm(\mu)$ converge to  the corresponding coefficients of the system (\ref{Dirac2}). As a consequence, choosing
$a_2(\mu) = (\sinh \alpha)^{-1}$,
the corresponding solution $a_1(\mu)$ and $a_3(\mu)$ of (\ref{Dirac10}) converge to $a_1$ and $a_3$ given by (\ref{Dirac3}).
Now the value of $\partial_y \psi_1(t,-1)$ is given by
(\ref{Dirac35}). 
Thus, by continuity, for any fixed $T > 0$, $\partial_y \psi_1(t,-1,\mu)$ converges uniformly as $\mu$ goes to $0$ to $\partial_y \psi_1(t,-1)$, namely to the case $\mu = 0$.

Note that this regularization is independent on $\nu$. It provides
an horizontal velocity $u_{1,h}(t,x,-1)$ arbitrarily close to
(\ref{Dirac100}).

%%%%%%%%%%%%%%%%%%%%

\subsection{Study of the boundary layer}

%%%%%%%%%%%%%%%%%%%%

We now turn to the study of the boundary layer $u_1^b(t,x,Y)$. We insert the Ansatz in Navier-Stokes equations and follow the classical derivation of Prandtl equations.  Let $Y$ be of order $O(1)$, and let us define $v_1$ and $v_2$ by 
$$ \begin{aligned}
%v(t,x,Y) &= (U_s(y),0) + \nu^\beta u_1(t,x,y) + \nu^\beta %u_1^b(t,x,Y)
%u^\nu &= (- 1 + \sqrt{\nu} Y,0) + \nu^\beta u_1(t,x,-1 +  %\sqrt{\nu} Y) + \nu^\beta u_1^b(t,x,Y)
u^\nu(t,x,y) & = (-1,0) 
+ \nu^\beta \Bigl( v_1(t,x,Y), \nu^{1/2} v_2(t,x,Y) \Bigr) 
+ \Bigl( O(\nu^{2 \beta}) + O(\sqrt{\nu}), O(\nu) \Bigr),
\end{aligned}
$$
since, as usual in boundary layers, using the divergence free condition, the vertical velocity is a factor  $\sqrt{\nu}$  smaller than the horizontal velocity. 
For $Y$ of order $O(1)$, the equation on the horizontal velocity gives, at order $O(\nu^\beta)$,
$$
\partial_t v_1 - \partial_x v_1 - \partial_Y^2 v_1 + \partial_x p = 0
$$
since $v_1 \partial_x v_1$ and $v_2 \partial_Y v_1$ are both negligible with respect to $\nu^\beta$.
The equation on the vertical velocity leads to
$$
\partial_Y p = 0,
$$
namely to $p = 0$, where $p$ is the corresponding pressure since $p$ vanishes at infinity.
Thus, as usual in boundary layer theory, the pressure does not change at leading order in the boundary layer.

The divergence free condition gives
$$
\partial_x v_1 + \partial_Y v_2 = 0.
$$
We are thus lead to solve
\beq \label{linPrandtl}
\partial_t u_{1,h}^b - \partial_x u_{1,h}^b - \partial_Y^2 u_{1,h}^b = 0
\eeq
with the boundary condition 
\beq \label{linPrandtl2}
u_{1,h}^b(t,x,0) = - u_{1,h}(t,x,-1).
\eeq
In order to remove the term $- \partial_x u_{1,h}^b$ in (\ref{linPrandtl}), we change the $x$ coordinates and from now on we work in a frame which moves with the flow at $y = -1$, namely with velocity $-1$. Navier-Stokes equations are invariant under this change of frame, and now the boundary condition is the usual Dirichlet condition $u = 0$ at $y = -1$.
Let $w(t,x,Y) = u_{1,h}^b(t,x-t,Y)$, then
$$
\partial_t w - \partial_Y^2 w = 0
$$
with in particular, when $\omega_1(0,y)$ is given by (\ref{Dirac}),
$$
w(t,0,0) = - \cos {\alpha t } 
-   a_1 \sinh {3 \alpha \over 2}  \cos {\alpha t \over 2}  -  a_3 \sinh {\alpha \over 2}   \cos  {3 \alpha t \over 2}   .
$$
Note that the boundary condition is periodic of period $4 \pi / \alpha$.
When $\mu$ goes to $0$, $w(t,0,0)$ converges to this explicit formula uniformly on every compact set in time.  Also note that, the boundary value of $v$ will increase then decrease as $t$ getting larger, namely, $\partial_Y^2w=\partial_t w$ at the boundary will change the sign which creates an inflection point for $Y>0$.

 The various profiles $u_j^b$ and $u_j^{b,1}$ satisfy linearized Prandtl equations with similar source terms. The construction of an asymptotic
expansion of the form (\ref{Ansatz}) is standard and we will not detail it here.
 We fix $N$ large enough in the Ansatz.
By construction, $u^\nu$ satisfies Navier-Stokes equation up to a small error term. We define $f^\nu$ to be this error term.

This ends the proof of the Theorem \ref{theorem1}, since, as soon as $t > 0$, $u_{1,h}^b$ does not vanish, hence $\nabla u^\nu$ is of order $\nu^{\beta - 1/2}$ is the boundary layer, which gives (\ref{ii4}) and (\ref{ii5}).

Note that $w(t,0,Y)$  can be explicitly computed using the Green function of the heat equation (double layer potentials ).
In particular, $w(t,0,Y)$ is rapidly decaying at infinity, like the Gaussian, and in particular faster than any exponential.
 It can also be easily numerically computed using for instance an implicit scheme.

Let us now discuss the evolution of the boundary layer in the particular case $\alpha = 1$ and in the limit case where $\omega_1(0,y)$ is given by (\ref{Dirac}). Let
$$
V_s(t,Y) =   u_{1,h}(t,-t,-1) + u_{1,h}^b(t,-t,Y). 
$$
We note that $u_{1,h}(t,-t,-1)$ oscillates in time. At first it decays with time. As a consequence, in small time,
$V_s(t,Y)$ is convex in $Y$.
However, later, the speed at infinity increases with time. This creates an inflexion point near the boundary and even an ``overshoot" of the flow.

%%%%%%%%%%%%%%%%%%%%%%%%%%%%%%%%%%%%%%%%%%%%%%

\section{Numerical investigation of the claim}

%%%%%

In this section, we numerically show  that for some value of $t_0$, the boundary layer at $x_0 = -t_0$ 
$$
V_s(t_0,Y) =  u_{1,h}(t_0,x_0,-1) + u_{1,h}^b(t_0,x_0,Y) 
$$
is unstable for linearized Euler equations. Note that $V_s(t_0,Y)$ has been computed in the previous section.

We look for an instability of the form
$$
v_{Ray}(t,x,Y) = \nabla^\perp \Bigl[ e^{i \alpha_{Ray} (x - c t) } \psi_{Ray}(Y) \Bigr].
$$
In the sequel, we choose $\alpha_{Ray}=\pm \sqrt{0.1}$.
We thus study the corresponding Rayleigh equation
\beq \label{Ray1}
(\nu^\beta V_s - c) (\partial_Y^2 - \alpha_{Ray}^2) \psi_{Ray} - \nu^\beta V_s'' \psi_{Ray} = 0.
\eeq 
 We rescale $c$ by
 $$
 c = \nu^\beta c_{Ray},
 $$
 which leads to the usual Rayleigh equation
 \beq \label{Ray2}
( V_s -  c_{Ray}) (\partial_Y^2 - \alpha_{Ray}^2) \psi_{Ray} - V_s'' \psi_{Ray} = 0,
\eeq 
\beq \label{Ray3}
\psi_{Ray}(0) = 0, \qquad \lim_{Y \to + \infty} \psi_{Ray}(Y) = 0.
\eeq
 We thus  look for a solution $(c_{Ray},\psi_{Ray})$ of (\ref{Ray2},\ref{Ray3}) with $\Im c_{Ray} > 0$, where $V_s$ is explicitly given.
The spectrum of Rayleigh operator is composed of two parts. First a continuous spectrum, which is the range of $V_s$, and lies on the real axis. There may also exist eigenvalues, which come in conjugate pairs. There may also exist embedded eigenvalues.
 
 Unfortunately, there are few theoretical results to study the existence of unstable eigenvalues, namely eigenvalues with $\Im c_{Ray}\neq 0$. Rayleigh criterium and its improvement by Fjortjoft, can not be applied since $V_s(t_0,Y)$ has an inflexion point, and there is no theoretical way to construct an unstable mode with a profile as complex as $V_s$. We thus have to rely on numerical computations.

 There are two ways to study numerically (\ref{Ray2},\ref{Ray3}). We can see it as a spectral problem or as a shooting problem. In order to obtain reliable numerical results, we follow both approaches.

To compute the spectrum of Rayleigh, we rewrite (\ref{Ray1}) under the form
$$
V_s \omega_{Ray} - V_s'' (\partial_Y^2 - \alpha^2)^{-1} \omega_{Ray} = c_{Ray} \omega_{Ray}
$$
where 
$$
\omega_{Ray} = (\partial_Y^2 - \alpha^2) \psi_{Ray}.
$$
We thus introduce the operator
$$
Ray \, \omega := V_s \omega - V_s'' A \omega
$$
where $A$ is the inverse of the Laplace operator, namely
$$
(\partial_Y^2 - \alpha_{Ray}^2) A \omega = \omega
$$
with  boundary conditions $A \omega(0) = 0$ and $A \omega(Y) \to 0$ as $Y \to + \infty$.
Then $Ray$, the Rayleigh operator in vorticity formulation, is a perturbation of the inverse of the Laplace operator.

We now have to  numerically evaluate the spectrum of $Ray$.
For this we choose some large $Y_0$ ($Y_0 = 30$ in our computations) and discretize $0 \le Y \le Y_0$ using a small step $h$ ($h = 0.001$). Let $\omega_k$ be the approximation of $\omega$ at $x_k = h k$. We approximate $A$ by the classical finite difference scheme
$$
\partial_Y^2 \omega(x_k) \approx {\omega_{k+1} + \omega_{k-1} - 2 \omega_k \over h^2} 
$$
for $1 \le k \le N$, together with the boundary conditions
$\omega_0 = 0$ (Dirichlet boundary condition at $y = 0$)
and $\omega_{N+1} = \omega_N$ (Neumann  boundary condition at $Y = N$). This approximation leads to a $N \times N$ matrix, whose spectrum can be numerically evaluated.
This gives a first numerical evaluation $\tilde c_{Ray}$ of the most unstable eigenvalue of $Ray$.

To check these numerical computations, we follow a completely different approach.
We  rewrite (\ref{Ray2}) as an ordinary differential equation 
\beq \label{Ray10}
\partial_Y^2 \psi = \alpha_{Ray}^2 \psi + {V_s''(t_0,Y) \over V_s(t_0,Y) - c} \psi 
\eeq
with $\psi(Y) \to 0$ as $Y \to + \infty$.
We note that $V_s''(Y)$ converges more than exponentially fast to $0$.
Classical results on ordinary differential equations then provide the existence of two solutions $\psi_\pm(Y,c)$, such that
$$
\psi_\pm(Y,c) \sim e^{\pm \alpha_{Ray} Y}
$$
as $Y \to + \infty$.
As we are looking for eigenmodes $\psi_{Ray}$ which go to $0$ at infinity, $\psi_{Ray}$ must be a multiple of $\psi_-$. Up to the multiplication by a constant, we may thus assume that $\psi_{Ray} = \psi_-(Y,c_{Ray})$.

Thus, $c$, with $\Im c > 0$ is an eigenvalue of Rayleigh equation if and only if
\beq \label{Ray11}
\psi_-(0,c) = 0.
\eeq
This is a shooting problem: we have to adjust $c$ such that (\ref{Ray11}) is satisfied.
We numerically solve (\ref{Ray10}) and then use a Newton method to determine its zero.

The numerical integration of  (\ref{Ray10}) is standard. For this we choose some large $Y_0$ ($Y_0 = 30$ in our computations), and solve (\ref{Ray10}) backwards to $Y = 0$, using a classical Runge-Kutta algorithm, starting from 
$\psi(Y_0,c) = e^{- \alpha_{Ray} Y_0}$ and $\partial_Y \psi(Y_0,c) = - \alpha_{Ray} e^{- \alpha_{Ray} Y_0}$. 
This procedure gives a numerical approximation of $\psi(0,c)$.
Note that this computation can be arbitrarily precise, provided the time step is small enough.

\begin{remark}
Fix $\alpha>0$. 
    Notice that here we focus on the case $c\notin \mathrm{Ran}\, V_s(t_0,Y)$. Thus the classical ODE argument together with the fact that $|V_s''(t_0,Y)|\lesssim e^{-c_1Y}$ for $c_1>0$ gives us that, there are two linearly independent solutions $\psi_1$ and $\psi_2$, such that at $Y_0$, 
    \begin{align*}
        &\psi_1(Y_0)= e^{-\alpha Y_0},\quad \partial_Y\psi_1(Y_0)=-\alpha e^{-\alpha Y_0},\\
        &\psi_2(Y_0)= e^{\alpha Y_0},\quad \partial_Y\psi_2(Y_0)=\alpha e^{\alpha Y_0}.
    \end{align*}
    Then 
    \begin{align*}
        \psi_1(Y)=e^{-\alpha Y}+e^{-\alpha Y}\int_{Y_0}^Y\int_{Y_0}^y \frac{V_s''(t_0,z)}{V_s(t_0, z)-c}\psi_1(z)e^{2\alpha y-\alpha z}dzdy,
    \end{align*}
    \begin{align*}
        \psi_2(Y)=e^{\alpha Y}+e^{\alpha Y}\int_{Y_0}^Y\int_{Y_0}^y \frac{V_s''(t_0,z)}{V_s(t_0, z)-c}\psi_2(z)e^{-2\alpha y+\alpha z}dzdy,
    \end{align*}
    which together with a classical fixed point argument gives that for $Y>Y_0$ with $Y_0$ large enough, 
    \begin{align*}
        |\psi_2(Y)-e^{\alpha Y}|\leq \delta_{Y_0} e^{\alpha Y}
    \end{align*}
    for some $0<\delta_{Y_0}<1$. Thus 
    \begin{align*}
        \psi_1(Y)=-2\alpha\psi_2(Y)\int_{Y_0}^{Y}\frac{1}{\psi_2(Z)^2}dZ+e^{-2\alpha Y_0}\psi_2(Y)
    \end{align*}
    and
    \begin{align*}
        \psi(Y)
        &=-2\alpha\psi_2(Y)\int_{\infty}^Y\frac{1}{\psi_2(z)^2}dz\\
        &=-2\alpha\psi_2(Y)\int_{Y_0}^Y\frac{1}{\psi_2(z)^2}dz-2\alpha \psi_2(Y)\int_{\infty}^{Y_0}\frac{1}{\psi_2(z)^2}dz.
    \end{align*}
    It is easy to check that, if $\psi(0)=0$ for some $c$, then $c$ is the eigenvalue and $\psi$ is the associated eigen-function. 
    We also have $|\psi_1(0)-\psi(0)|\to 0$ as $Y_0\to \infty$. In the numerical computation, we find zeros for $\psi_1(0)$ instead of $\psi(0)$ and take $Y_0$ large enough. 
\end{remark}

We now use a Newton algorithm to solve (\ref{Ray11}).
We start from $\tilde c_{Ray}$, given through the spectral approach, which is an approximate eigenvalue of (\ref{Ray11}).

Newton algorithm gives another numerical approximation $\hat c_{Ray}$. It turns out that $\tilde c_{Ray}$ and $\hat c_{Ray}$ are very close, which validates the numerical computations, since the same result has been obtained by two different methods sharing no common computer lines code.

\medskip

\begin{figure}[h!]
\centerline{\includegraphics[width=11cm]{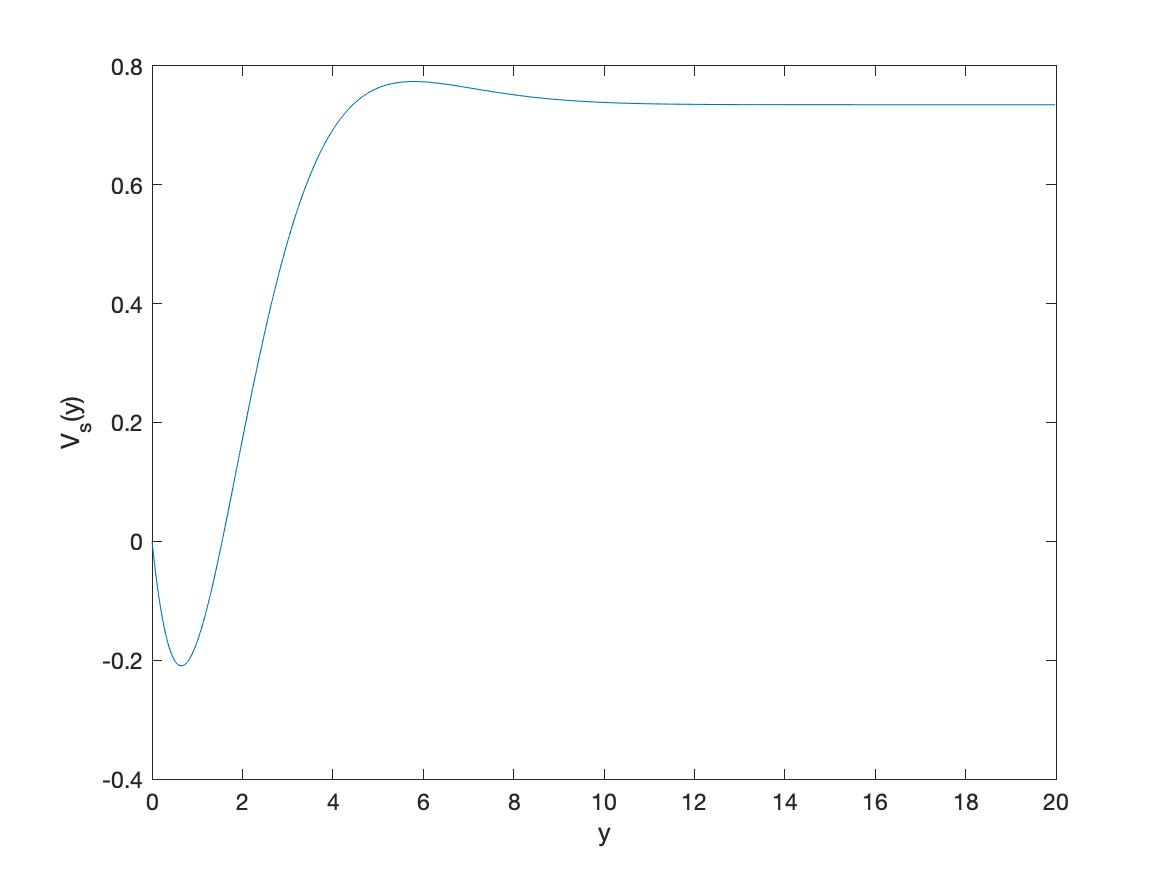}}
\caption{The boundary layer velocity $V_s(t_0,Y)$ at $t_0 = 7.65$}
\label{Us}
\end{figure}

Let us now detail the numerical computations
 For $\alpha = 1$, at $t_0 = 7.65$, the velocity profile 
is displayed on figure \ref{Us}. This profile is not convex and has an inflexion point.

\begin{figure}[h!]
\centerline{\includegraphics[width=11cm]{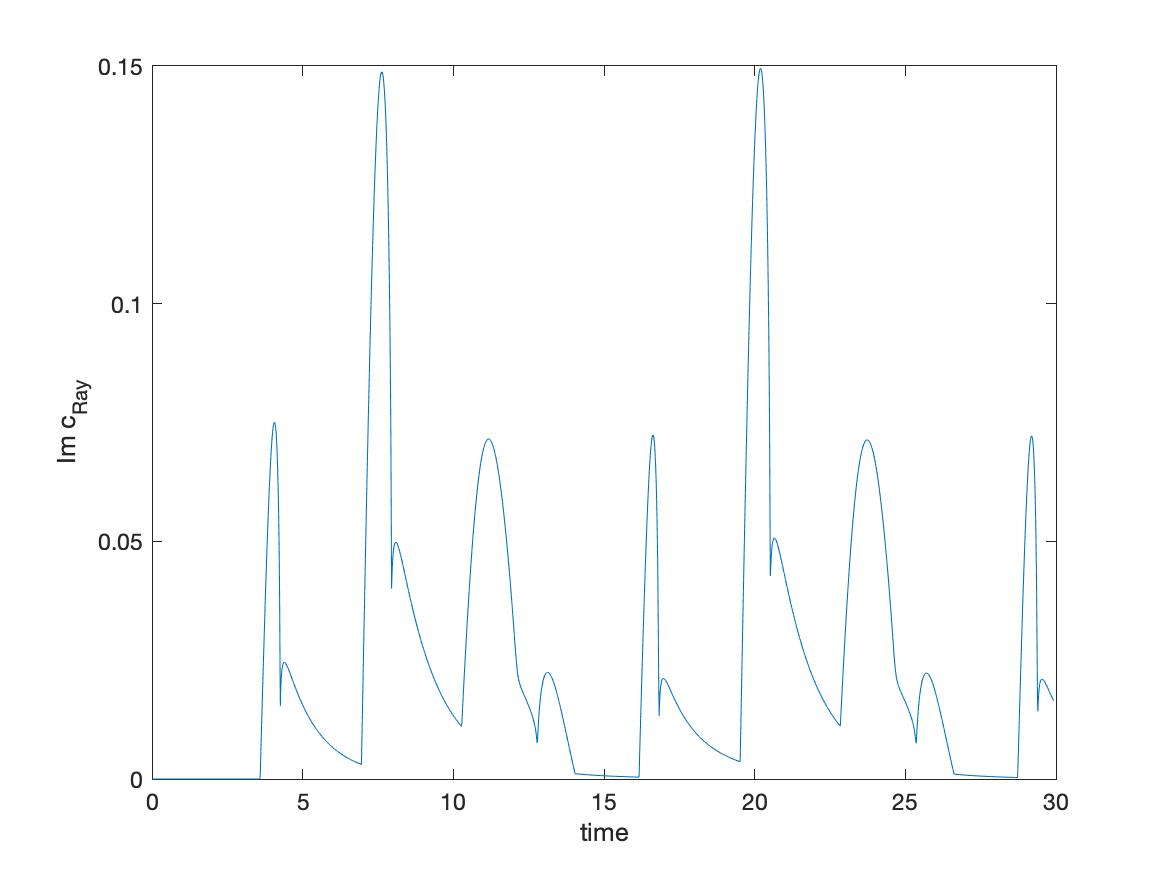}}
\caption{Speed of growth of the most unstable mode as a function of time}
\label{evolutionunstable}
\end{figure}

\begin{figure}[h!]
\centerline{\includegraphics[width=10cm]{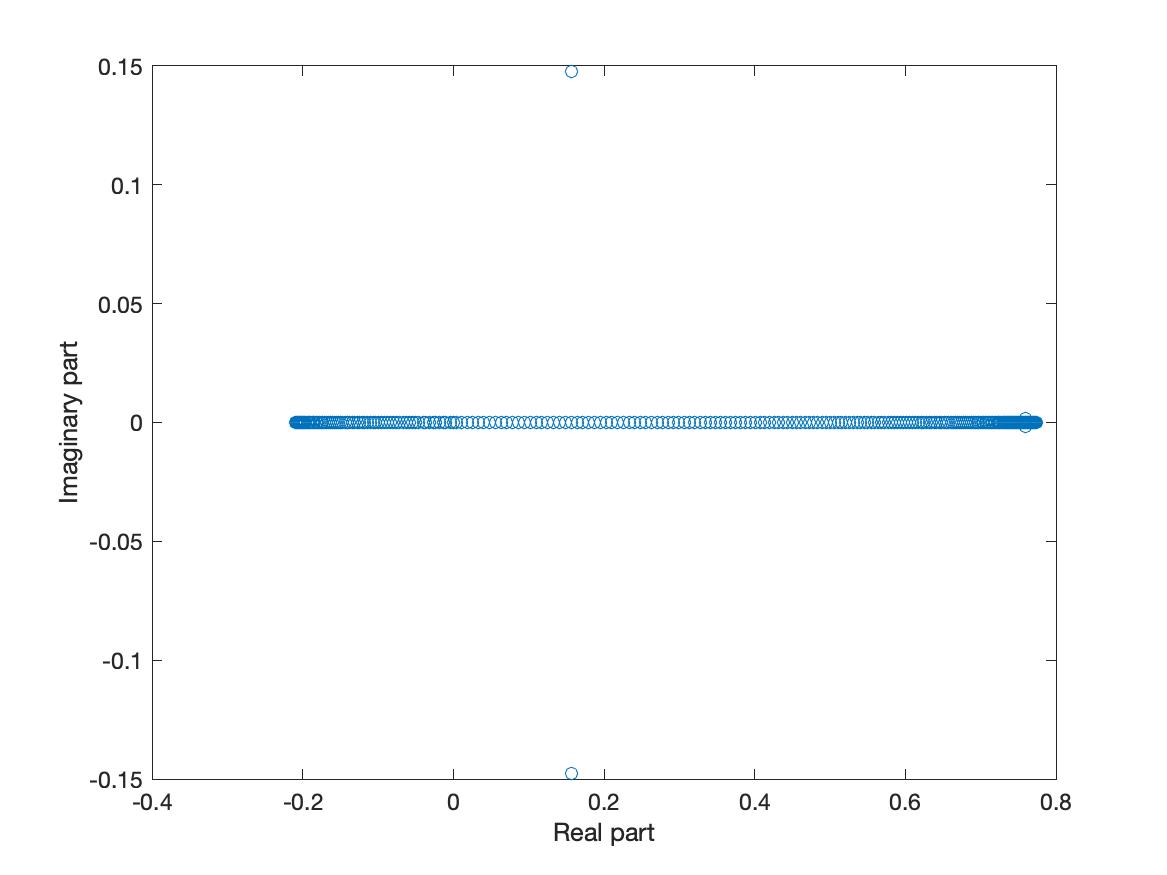}}
\caption{Numerical approximation of the spectrum of Rayleigh at $t_0 = 7.65$}
\label{spectrum}
\end{figure}

\begin{figure}[h!]
\centerline{\includegraphics[width=10cm]{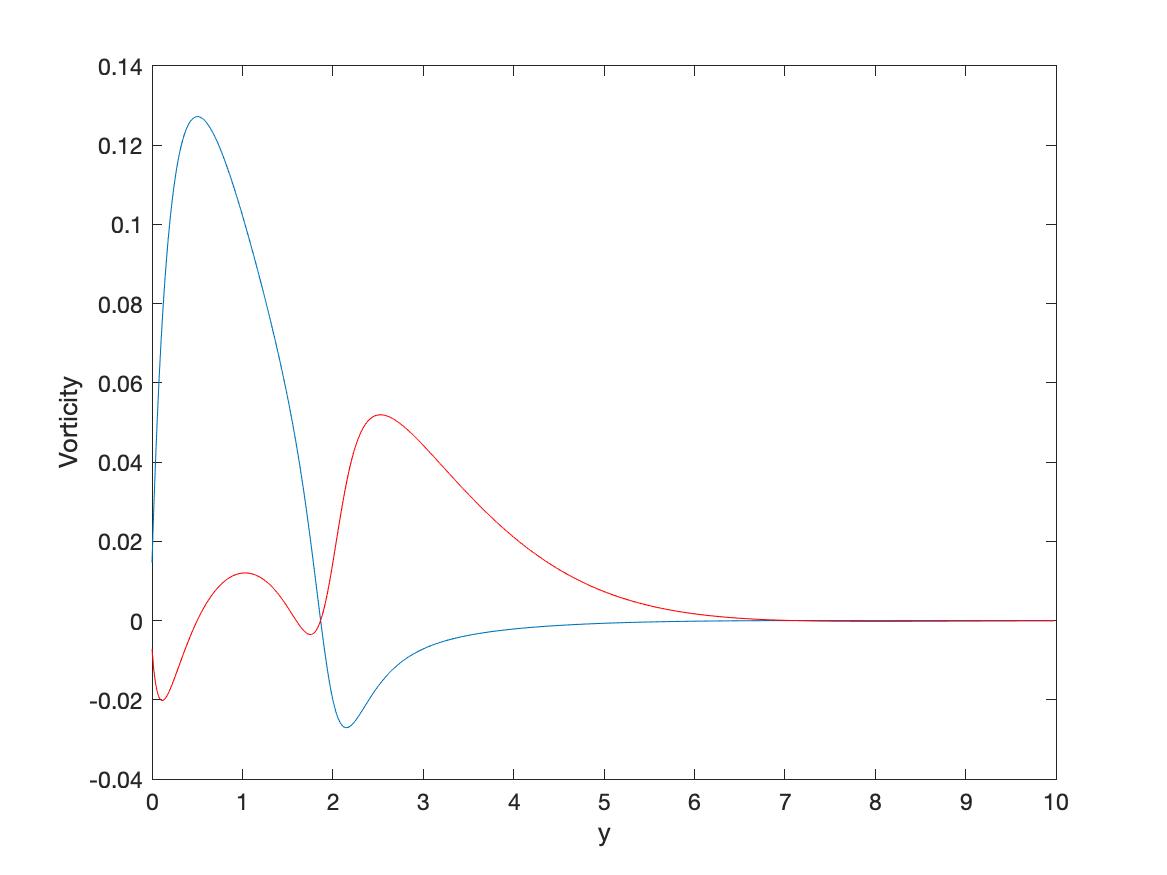}}
\caption{Vorticity $\omega_{Ray}$ of the eigenfunction at $t_0 = 7.65$ (real part in blue, imaginary part in red)}
\label{vorticity}
\end{figure}

\begin{figure}[h!]
\centerline{\includegraphics[width=10cm]{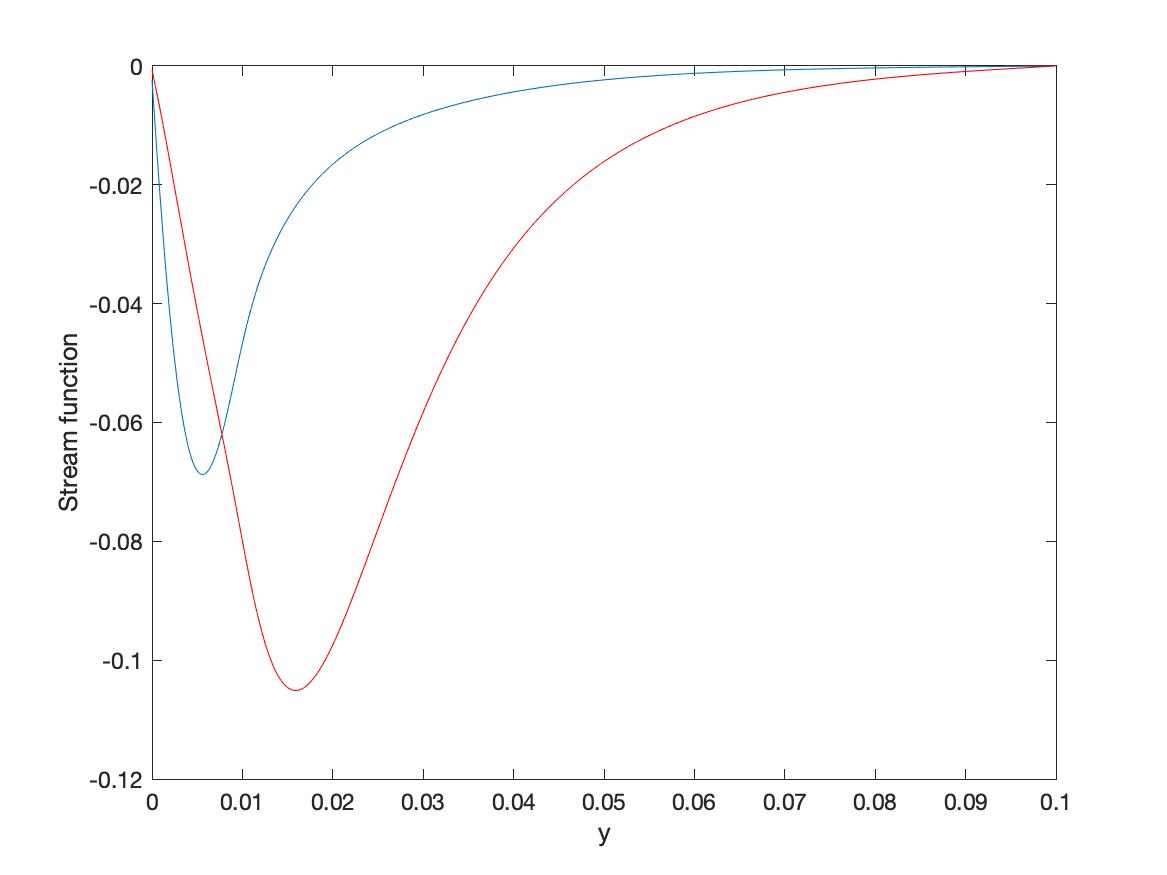}}
\caption{Unstable eigenfunction  $\psi_{Ray}$ at $t_0 = 7.65$ (real part in blue, imaginary part in red)}
\label{eigenfunction}
\end{figure}

We choose $\alpha_{Ray} = \sqrt{0.1}$.
Figure \ref{evolutionunstable} diplays the imaginary part of the most unstable eigenvalue of Rayleigh operator, if there is one, and $0$ if all the eigenvalues are real.
An unstable mode appears at $t \approx 4$. The rate of growth of the instability is maximal at $t_0 \approx 7.6$ and then decays before growing again.
Figure \ref{spectrum} shows the spectrum of the corresponding Rayleigh equation. 
There exist two conjugate eigenvalues which are not real. The continuous spectrum is the range of $U_s$, between approximately $-0.2$ and $0.8$.
The real and imaginary parts of the unstable vorticity $\omega_{Ray}$ and of the eigenfunction $\psi_{Ray}$ are displayed on figures \ref{vorticity} and \ref{eigenfunction}.

The output of this numerical analysis is a solution $(c_{Ray},\psi_{Ray})$ of the Rayleigh equations (\ref{Ray2},\ref{Ray3}) such that
$\Im c_{Ray} > 0$.
Figure \ref{eigenfunction} shows the eigenfunction $\psi_{Ray}(Y)$. Note that $\psi_{Ray}$ goes exponentially fast to $0$ at infinity, and that $\psi_{Ray}(0) = 0$.
We note that 
\beq \label{sub}
\partial_Y \psi_{Ray}(0) \ne 0.
\eeq
The unstable eigenvalues of Rayleigh operator depend smoothly on the profile $V_s$. Thus, provided $\mu$ is small enough, the profile $V_s$ has an unstable eigenvalue which is close to $c_{Ray}$ with a corresponding eigenmode which is close to $\psi_{Ray}$.

%%%%%

\section{Proof of Theorem \ref{theomain} }

%%%%%%%%%%%%%%

\subsection{Viscous linear instability in the boundary layer}

%%%%%%%%%%%%%%

Using the claim, we obtain that
$$
V_s(t_0,Y) =  u_{1,h}(t_0,x_0,-1) + u_{1,h}^b(t_0,x_0,Y) 
$$
is linearly unstable for Euler equation, with corresponding exponentially growing solution
$$
v_{Ray}(T,X,Y) = \nabla^\perp \Bigl[ e^{i \alpha_{Ray} (X - \nu^\beta c_{Ray} T) } \psi_{Ray}(Y) \Bigr].
$$
The next step is to construct an approximate unstable mode $(c_{Orr},\psi_{Orr})$ for the Orr-Sommerfeld equations
which are, taking into account the $\nu^\beta$ factor,
$$
(V_s - c) (\partial_Y^2 - \alpha_{Ray}^2) \psi - V_s'' \psi = {i \tilde \nu \over \alpha_{Ray} \nu^\beta} (\partial_Y^2 - \alpha_{Ray}^2)^2 \psi.
$$ 
Let
$$
\hat \nu = {\tilde \nu \over \nu^\beta} = \nu^{1/2 - \beta}.
$$
As $\partial_Y \psi_{Ray}(0) \ne 0$, we have to add a new boundary layer in order to recover the boundary condition
$\psi_{Orr}(0) = \partial_Y \psi_{Orr}(0) = 0$.
The construction is classical \cite{GN2}. We look for $\psi_{Orr}$ and $c_{Orr}$ under the form
\beq \label{porr1}
\psi_{Orr}(Y) = \sum_{n \ge 0} \hat \nu^{n/2} \psi_n(Y) + \hat \nu^{n/2} \psi_n^b(\tilde{Y}) 
\eeq
where
$$
\tilde{Y} = {Y \over \hat \nu^{1/2}} = {Y \over \nu^{1/4 - \beta/2}} = {y + 1 \over \nu^{3/4 - \beta/2}} .
$$
We also look for $c_{Orr}$ under the form
\beq \label{porr2}
c_{Orr} = \sum_{n \ge 0} \hat \nu^{n/2} c_n.
\eeq
The sums (\ref{porr1}) and (\ref{porr2}) are just asymptotic expansions. In the sequel, we truncate them at some large integer $M$, which leads
to approximate eigenvectors and eigenvalues of Orr-Sommerfeld equation.

We start with $\psi_0 = \psi_{Ray}$, $\psi_0^b = 0$ and $c_0 = c_{Ray}$.
The construction of the various $\psi_n$, $\psi_n^b$ and $c_n$ is straightforward. We insert these Ansatz in Orr-Sommerfeld equation,
which gives a series of equation, one by power of $\hat \nu^{1/2}$. This provides, after truncation, an approximate solution of Orr-Sommerfeld equation.
We do not detail the construction any further.

%%%%%%%%%%%%%%

\subsection{Nonlinear instability in the boundary layer}

%%%%%%%%%%%%%%

We now turn to the construction of the nonlinear instability of the boundary layer.
We fix $t_0 = 7.65$ and $x_0 = - t_0$, and construct a linear instability at this time and this point.
First we make an isotropic change of variables
$$
(T,X,Y) = {t - t_0,x - x_0,y  + 1 \over \sqrt{\nu}} .
$$
Navier Stokes equations remain invariant, expect for the viscosity which is now 
$$
\tilde \nu = \nu^{1/2}.
$$
In these new variables $u^\nu$ becomes
$$
\tilde u^\nu(T,X,Y) = u^\nu(t_0 + \nu^{1/2} T, x_0 + \nu^{1/2} X, - 1 + \nu^{1/2} Y) 
$$
and slowly depends on the rescaled variables $T$, $X$ and $Y$.

We now construct an instability in the boundary layer by following the general strategy of \cite{BGlayers}.
In the $(T,X,Y)$ variable, $\tilde u^\nu$ is changing in times of order $\nu^{-1/2}$, in $X$ of order $\nu^{-1/2}$ and in $Y$ of orders $1$ (boundary layer) and
 $\nu^{-1/2}$ (interior behavior).

In the previous paragraph, we have constructed an unstable mode for linearized Navier-Stokes equations around $V_s(Y)$. 
This exponentially growing solution is of the form
$$
v_1(T,X,Y) = \nabla^{\perp} \Bigl[ e^{i \alpha_{Ray} (X - \nu^\beta c_{Orr} T)} \psi_{Orr}(Y) \Bigr].
$$
It grows over times $T$ of order $\nu^{-\beta}$ (corresponding to times $t$ of order $\nu^{1/2 - \beta}$), and is periodic in $X$ with a period 
$2 \pi / \alpha_{Ray}$ (corresponding to a very small period $2 \pi \nu^{1/2} / \alpha_{Ray}$ in the $x$ variables). Moreover
it has two scales in $Y$, namely $1$ and $\hat \nu^{1/2} = \nu^{1/4-\beta/2}$, 
corresponding to $\nu^{1/2}$ and $\nu^{3/4 -\beta/2}$ in the original $y$ variable.

Following Theorem $4.1$ of \cite{BGlayers} which is recalled in the Appendix, we can then construct a nonlinear instability starting form $v_1$.
This instability then almost reaches a size $\nu^{\beta}$.
This ends the proof of Theorem \ref{theomain}.

%%%%%%%%%%%%%%%%%%%%%%%%%%%%%%%%%%%%%%%%%

\subsection{Instability in $\Omega_\rit$  \label{boundary2}}

%%%%%%%%%%%%%%%%%%%%%%%%%%%%%%%%%%%%%%%%%

The proof in the whole strip $\Omega_\rit$ is obtained by "localizing" the previous proof. Namely, we look for $u_1(t,x,y)$ which is of the form
\beq \label{u11}
u_1(t,x,y) = \nabla^\perp \int \chi(\alpha) \psi_\alpha(t,y) e^{i \alpha x} d\alpha + \cc
\eeq
where $\chi(\alpha)$ is a smooth positive function, supported on $[1 - \eta,1 + \eta]$ for some small $\eta$, with unit integral.
The functions $\psi_\alpha$ are constructed as the function $\psi_1$.
Using a stationary phase theorem, we see that $u_1$ rapidly decays at infinity.
The end of the proof is similar.

%%%%%%%%%%%%%%%%%%%%%%%%%%%%%%

\section{Proof of Theorem \ref{theo3}}

%%%%%%%%%%%%%%%%%%%%%%%%%%%%%%

We now turn to the proof of Theorem \ref{theo3}. Let $u^\nu_0$ be a sequence of initial data such that
$$
\| u^\nu_0-(y,0) \|_{H^s} \le \nu^{\beta}
$$
with $\beta > 1/2$. The first step is to construct an approximate solution $u^{app}$ of the form (\ref{Ansatz})
on some time interval $[0,T]$. Note that $T$ may be arbitrarily large.
By construction, $u^{app}$ satisfies the Navier-Stokes equation, up to a very small error term
$$
\partial_t u^{app} + (u^{app} \cdot \nabla ) u^{app}
- \nu \Delta u^{app} + \nabla p^{app} = \nu^{N'} R^{app},
$$
$$
\nabla \cdot u^{app} = 0,
$$
where $R^{app}$ and $\partial_t R^{app}$ are uniformly bounded in $L^\infty([0,T],H^3)$.

%{\color{red} How about $\partial_tR^{app}$? }

Let $u^\nu$ be the genuine solution of Navier-Stokes equations with initial data $u^\nu_0$. Then
$$
v^\nu = u^\nu - u^{app}
$$
satisfies
$$
 \partial_t v^\nu + (u^{app} \cdot \nabla) v^{\nu}
+ (v^\nu \cdot \nabla) u^{app} + (v^\nu \cdot \nabla) v^\nu- \nu \Delta v^\nu + \nabla p^\nu = - \nu^{N'} R^{app}+f^{\nu},
$$
$$
\nabla \cdot v^\nu = 0.
$$
Moreover, by construction of the Ansatz, $v^\nu = 0$ at $t = 0$.
We note that 
$$
\| \nabla u^{app} \|_{L^\infty} \lesssim 1+\nu^{\beta - 1/2} \le C
$$
thus, a classical energy estimate gives
$$
{1 \over 2} \partial_t \| v^\nu \|_{L^2}^2 
+ \nu \| \nabla v^\nu \|_{L^2}^2 \le C \| v^{\nu} \|_{L^2}^2 
+ C \nu^{2N'}
$$
and, for any $0 \le t \le T$, using Gronwall inequality,
$$
\| v^\nu(t) \|_{L^2}^2 + \nu \int_0^t \| \nabla v^\nu(\tau) \|_{L^2}^2 d\tau \le C \nu^{2 N'}.
$$
We also have $\|\langle \partial_x \rangle^2 \partial_x u^{app}\|_{L^{\infty}}\lesssim \nu^{\beta}$ and $\|\langle \partial_x \rangle^3\nabla u^{app}\|_{L^{\infty}}\lesssim 1$, where 
\begin{align*}
    \|\langle \partial_x\rangle^m f\|_{L^2}:=\left(\sum_{k\in \mathbb{Z}}\|(1+\alpha^2)^{\frac{m}{2}}\widehat{f}(\alpha,\cdot)\|_{L^2_y}^2\right)^{\frac{1}{2}}\approx \sum_{i=0}^m\|\partial_x^i f\|_{L^2}
\end{align*}
and $\widehat{f}$ is the Fourier transform of $f$ in $x$ variable. 
Thus by a similar argument, we obtain that
\begin{align*}
    \frac{1}{2}\frac{d}{dt}\|\partial_x v^{\nu}\|_{L^2}^2+\nu\|\nabla \partial_x v^{\nu}\|_{L^2}^2\lesssim &\|\partial_x v^{\nu}\|_{L^2}^2
    +\nu^{2\beta}\|\nabla v^{\nu}\|_{L^2}^2\\
    &+\|\partial_xv^{\nu}\|_{L^2}\|\nabla \partial_xv^{\nu}\|_{L^2}\|\nabla v^{\nu}\|_{L^2}+\nu^{2N'},\\
    \frac{1}{2}\frac{d}{dt}\|\partial_{xx} v^{\nu}\|_{L^2}^2+\nu\|\nabla \partial_{xx} v^{\nu}\|_{L^2}^2\lesssim &\|\langle \partial_{xx}\rangle v^{\nu}\|_{L^2}^2
    +\nu^{2\beta}\|\nabla v^{\nu}\|_{L^2}^2+\nu^{2\beta}\|\nabla \partial_xv^{\nu}\|_{L^2}^2\\
    &+\|\partial_xv^{\nu}\|_{L^2}^{\frac{1}{2}}\|\nabla \partial_xv^{\nu}\|_{L^2}^{\frac{3}{2}}\|v_{xx}^{\nu}\|_{L^2}^{\frac{1}{2}}\|\nabla \partial_x^2 v^{\nu}\|_{L^2}^{\frac{1}{2}}\\
&+\|\partial_{x}^2v^{\nu}\|_{L^2}\|\nabla \partial_x^2v^{\nu}\|_{L^2}\|\nabla  v^{\nu}\|_{L^2}+\nu^{2N'},
    \\
 \frac{1}{2}\frac{d}{dt}\|\partial_{x}^3 v^{\nu}\|_{L^2}^2+\nu\|\nabla \partial_{x}^3 v^{\nu}\|_{L^2}^2\lesssim &\|\langle \partial_{x}\rangle^3 v^{\nu}\|_{L^2}^2
    +\nu^{2\beta}\| \langle \partial_x\rangle^2\nabla v^{\nu}\|_{L^2}^2+\nu^{2\beta}\|\nabla v^{\nu}\|_{L^2}^2\\
&+\|\langle\partial_x\rangle^2v^{\nu}\|_{L^2}^{\frac{1}{2}}\|\nabla \langle\partial_x\rangle^2v^{\nu}\|_{L^2}^{\frac{3}{2}}\|v_{xxx}^{\nu}\|_{L^2}^{\frac{1}{2}}\|\nabla \partial_x^3 v^{\nu}\|_{L^2}^{\frac{1}{2}}    
\\&
+\|\partial_{x}^3v^{\nu}\|_{L^2}\|\nabla \partial_x^3v^{\nu}\|_{L^2}\|\nabla  v^{\nu}\|_{L^2}+\nu^{2N'}
\end{align*}
which gives that 
\begin{align*}
    \|\langle \partial_x\rangle^3 v^{\nu}\|_{L^{\infty}L^2}+\nu \|\langle \partial_x\rangle^3 \nabla v^{\nu}\|_{L^2L^2}\lesssim \nu^{N'}. 
\end{align*}
We now estimate $\partial_t v^\nu$. We have
$$
\partial_t \partial_t v^\nu + (u^{app} \cdot \nabla) \partial_t v^\nu + (\partial_t u^{app} \cdot \nabla) v^\nu
+ (\partial_t v^\nu \cdot \nabla) u^{app} 
+ (v^\nu \cdot \nabla) \partial_t u^{app} 
$$
$$
+ (\partial_t v^\nu \cdot \nabla) v^\nu
+ (v^\nu \cdot \nabla) \partial_t v^\nu
- \nu \Delta \partial_t v^\nu + \nabla \partial_t p^\nu= 
- \nu^{N'} \partial_t R^{app}.
$$
We fulfill $L^2$ energy estimates on $\partial_t v^\nu$. We have, using that $\| \nabla u^{app} \|_{L^\infty}$ is bounded,
$$
{1 \over 2} \partial_t \| \partial_t v^\nu \|_{L^2}^2 
+ \nu \| \nabla \partial_t v^\nu \|_{L^2}^2 \le
\int | \partial_t v^\nu (\partial_t u^{app} \cdot \nabla) v^\nu |
+ C \| \partial_t v^\nu \|_{L^2}^2 
+ \int | \partial_t v^\nu v^\nu \cdot \partial_t \nabla u^{app} |
$$
$$
+ \int | \partial_t v^\nu |^2 |\nabla v^\nu |
+ \nu^{2N'} \| \partial_t R^{app} \|_{L^2}^2.
$$
We have, using that $\partial_t u^{app}$ is uniformly bounded in $\nu$,
$$
\Bigl| \int \partial_t v^\nu (\partial_t u^{app} \cdot \nabla) v^\nu \Bigr| \le C \| \partial_t v \|_{L^2} \| \nabla v^\nu \|_{L^2} 
\le C \| \partial_t v^\nu \|_{L^2}^2 + C \| \nabla v^\nu \|_{L^2}^2
$$
where the second term is bounded by $\nu^{N'-1}$ in $L^2([0,T])$ norm.
Moreover, as $\nu^{1/2} \| \partial_t \nabla u^{app} \|_{L^\infty}$ is uniformly bounded,
$$
\Bigl| \int \partial_t v^\nu v^\nu \cdot \partial_t \nabla u^{app} \Bigr| \le C \nu^{-1/2} \| \partial_t v^\nu \|_{L^2} \| v^\nu \|_{L^2} \le C \| \partial_t v^\nu \|_{L^2}^2 + C \nu^{-1} \| v^\nu \|_{L^2}^2,
$$
where the second term is bounded by $C \nu^{N-1}$ in $L^\infty([0,T])$.
Next
$$
\int | \partial_t v^\nu|^2 | \nabla v^\nu |
\le \| \partial_t v^\nu \|_{L^4}^2 \| \nabla v^\nu \|_{L^2}
\le \| \partial_t v^\nu \|_{L^2} \| \partial_t \nabla v^\nu \|_{L^2} \| \nabla v^\nu \|_{L^2}
$$
$$
\le {\nu \over 2} \| \nabla \partial_t v^\nu \|_{L^2}^2 + C \nu^{-1} \| \nabla v^\nu \|_{L^2}^2 \| \partial_t v^\nu \|_{L^2},
$$
where $\nu^{-1} \| \nabla v^\nu \|_{L^2}^2$ is bounded in $L^1([0,T])$ by $\nu^{N-1}$.
Combining all these estimates, we obtain that
$$
\partial_t \|  \partial_tv^\nu \|_{L^2}^2
+ \nu \| \nabla\partial_t v^\nu \|_{L^2}^2
\le \Bigl( C + \nu^{N-1} \phi_1(t) \Bigr) \| \partial_t v^\nu \|_{L^2}^2 + \nu^{N-2} \phi_2(t)
$$
where $\phi_1$ and $\phi_2$ are bounded in $L^1([0,T])$.
Using Gronwall inequality and the fact that $\partial_tv|_{t=0}=-\nu^{N'}R^{app}+f^{\nu}$, we obtain that
$\partial_t v^\nu$ is bounded by $\nu^{N'-2}+\nu^{N}$ in $L^\infty([0,T],L^2) \cap L^2([0,T],H^1)$.

Now we use the Stokes estimates, namely, 
\begin{align}\label{eq:stokes-est}
    \nu\|\Delta v^{\nu}\|_{L^2}+\|\nabla p\|_{L^2}\lesssim \|F^{\nu}\|_{L^2}
\end{align}
where 
$$
F^\nu = -\partial_t v^{\nu}- (u^{app} \cdot \nabla) v^\nu - (v^\nu \cdot \nabla) u^{app}
- (v^\nu \cdot \nabla) v^\nu - \nu^{N'} R^{app}+f^{\nu}.
$$
Therefore, it holds that
\begin{align*}
    \|F^{\nu}\|_{L^2}\lesssim &\|\partial_t v^{\nu}\|_{L^2}+\nu^{\beta}\|\nabla v^{\nu}\|_
    {L^2}+\|v^{\nu}\|_{L^2}+\|v^{\nu}\|_{L^2}^{\frac{1}{2}}\|\nabla v^{\nu}\|_{L^2}\|\nabla^2 v\|_{L^2}^{\frac{1}{2}}+\nu^{N'}+\nu^{N}\\
    &\lesssim \nu^{N'-2}+ \nu^{\beta}\| v^{\nu}\|_
    {L^2}^{\frac{1}{2}}\|\nabla^2 v^{\nu}\|_{L^2}^{\frac{1}{2}}+\|v^{\nu}\|_{L^2}\|\nabla^2 v\|_{L^2}+\nu^{N'}+\nu^{N}
\end{align*}
which together with \eqref{eq:stokes-est} gives that
\begin{align*}
    \nu\|\Delta v^{\nu}\|_{L^2}+\|\nabla p\|_{L^2}\lesssim \nu^{N'-2}. 
\end{align*}
By taking the Fourier transform in $x$, we have
\begin{align*}
    \|\partial_y v^{\nu}\|_{L^{\infty}_{x,y}}
    &\lesssim \sum_\alpha\|\widehat{\partial_y v^{\nu}}(\alpha,\cdot)\|_{L^{\infty}_y}\\
    &\lesssim \sum_\alpha \|\langle \alpha\rangle^3\widehat{ v^{\nu}}(\alpha,\cdot)\|_{L^{2}_y}^{\frac{1}{4}}\|\widehat{ \partial_y^2 v^{\nu}}(\alpha,\cdot)\|_{L^{2}_y}^{\frac{3}{4}}\langle \alpha\rangle^{-\frac{3}{4}}\\
    &\lesssim \|\langle\partial_x\rangle^3 v^{\nu}\|_{L^2}^{\frac{1}{4}}
    \|\nabla^2 v^{\nu}\|_{L^2}^{\frac{3}{4}}\lesssim \nu^{N'-\frac{9}{4}}, 
\end{align*}
and
\begin{align*}
    \|\partial_x v^{\nu}\|_{L^{\infty}_{x,y}}
    &\lesssim \sum_\alpha \|\langle \alpha\rangle^3\widehat{v^{\nu}}(\alpha,\cdot)\|_{L^{2}_y}^{\frac{3}{4}}\|\widehat{ \partial_y^2 v^{\nu}}(\alpha,\cdot)\|_{L^{2}_y}^{\frac{1}{4}}\langle \alpha\rangle^{-\frac{5}{4}}\\
    &\lesssim \|\langle\partial_x\rangle^3 v^{\nu}\|_{L^2}^{\frac{3}{4}}
    \|\nabla^2 v^{\nu}\|_{L^2}^{\frac{1}{4}}\lesssim \nu^{N'-\frac{3}{4}}, 
\end{align*}
which gives $\|\nabla v^{\nu}\|_{L^{\infty}}\lesssim \nu^{N'-\frac{9}{4}}$.
This ends the proof of Theorem \ref{theo3}.

\section{The case of the pipe flow}

%%%%%%%%%%%%%%%%%%%%%

The pipe flow is the flow $U_s(r) = 1 - r^2$ in the cylinder of unit radius. This flow is physically known to be spectrally stable at any Reynolds number, exactly like the Couette flow.

Under axisymmetric perturbations without swirl, the vorticity equation is
\beq \label{vorticityPipe}
\partial_t \omega + (1 - r^2) \partial_z \omega -\nu(\partial_r^2+\frac{1}{r}\partial_r-\frac{1}{r^2}+\partial_z^2)\omega +u^r\partial_r\omega+u^z\partial_z\omega = 0,
\eeq
which has a similar linearized equation to \eqref{transportomega} when $\nu=0$. 
The computations are then similar.

%%%%%%%%%%%%%%%%%%%%%

\section{Appendix}

%%%%%%%%%%%%%%%%%%%%%

We now recall the Theorem $4.1$ of \cite{BGlayers}.
Let $u^\eps(t,x,y)$ be a sequence of solutions of Navier-Stokes equations which slowly depend on time,
 have "large structures" in the $x$ and $y$ variables, and a "boundary layer behavior" near $y = 0$,
 namely let $u^\eps(t,x,y)$ be a sequence of solutions of the form
\beq \label{ueps}
u^\eps(t,x,y) = 
\left( \begin{array}{c} U_1^{int,\eps}(\eps t,\eps x,\eps y) + U_1^{bl,\eps}(\eps t, \eps x, y) \cr
U_2^{int,\eps}(\eps t,\eps x,\eps y) + \eps U_2^{bl,\eps}(\eps t, \eps x, y)
\end{array} \right)
\eeq
where $\eps$, which depends on $\nu$, goes to $0$ as $\nu \to 0$. 
We assume that $U_1^{int,\eps}$ and $U_2^{int,\eps}$ are smooth in 
$T = \eps t$, $X = \eps x$ and $Y = \eps y$,
and converge to some functions $U_1^{int,0}$ and $U_2^{int,0}$ in $C^\infty$ as $\eps \to 0$.

 We moreover assume that $U_1^{bl,\eps}$ and $U_2^{bl,\eps}$ are smooth in  $T$, $X$ and $y$, 
 are exponentially decaying in $y$ (as well as all their derivatives),  and
 converge in $C^\infty$ (with uniform exponential decay) to some functions $U_1^{b,0}$ and $U_2^{b,0}$.
We define
$$
U_1^\eps(t,x,y) = U_1^{int,\eps}(\eps t,\eps x,\eps y) + U_1^{bl,\eps}(\eps t, \eps x, y)
$$
and
$$
U_2^\eps(t,x,y) = U_2^{int,\eps}(\eps t,\eps x,\eps y) + \eps U_2^{bl,\eps}(\eps t, \eps x, y),
$$
and similarly for $U_1^0$ and $U_2^0$.

\begin{theorem} \label{theogeneralfast}
Let us assume that $(U_1^{int,0}(0,0,0) + U_1^{b,0}(0,0,y),0)$ is a spectrally unstable shear layer for Euler equations, 
namely that there exists a
solution $(\alpha_{Ray},c_{Ray},\psi_{Ray})$ to the Rayleigh equation with $\Im c_{Ray} > 0$,
$\alpha_{Ray} \ne 0$,  $\psi_{Ray}(y) \not\equiv 0$ and $\psi_{Ray}(0) = 0$.
 
Assume that $\eps(\nu) \lesssim \nu^\gamma$ for some $\gamma > 0$. Then, for any arbitrarily small  positive $\theta$,  
for any arbitrarily large $N$ and arbitrarily large $s$, there exists a solution $v^\nu$ of Navier-Stokes equations with forcing term $f^\nu$,
 and a time $T^\nu$, such that, for $\nu$ small enough,
\beq \label{theogeneral1-2}
\| v^\nu(0,\cdot,\cdot) - u^{\eps(\nu)}(0,\cdot,\cdot) \|_{H^s} \le \nu^N,
\eeq
\beq \label{theogeneral2-2}
\| f^\nu \|_{L^\infty([0,T^\nu],H^s)} \le \nu^N
\eeq
and
\beq \label{theogeneral3-2}
\| u^\nu(T^\nu,\cdot,\cdot) - u^{\eps(\nu)}(T^\nu,\cdot,\cdot) \|_{L^\infty} \ge  \nu^\theta ,
\eeq
with 
$$
T^\nu \sim C_0 \log \nu^{-1}
$$
for some positive constant $C_0$.
\end{theorem}

\section*{Acknowledgements}
The work of N. Masmoudi is supported by NSF grant DMS-1716466 and by Tamkeen under the NYU Abu Dhabi Research Institute grant of the center SITE.
 D. Bian is supported by NSF of China under the grant 12271032.

%%%%%%%%%%%%%%%%%%%%%%%%%%
\bibliographystyle{abbrv.bst} 
\bibliography{references.bib}

%%%%%%%%%%%%%%%%%%%%%%%%%%%

\end{document}